\documentclass[12pt,leqno]{article}

\def\re{\mathop{\rm Re}}
\def\im{\mathop{\rm Im}}

\newtheorem{theorem}{Theorem}
\newtheorem{lemma}[theorem]{Lemma}
\newtheorem{proposition}[theorem]{Proposition}
\newtheorem{sublemma}[theorem]{Sublemma}
\newtheorem{definition}[theorem]{Definition}
\newtheorem{corollary}[theorem]{Corollary}
\newtheorem{problem}[theorem]{Problem}
\newtheorem{remark}[theorem]{Remark}
\newtheorem{claim}[theorem]{Claim}
\newtheorem{assumptions}[theorem]{Assumptions}
\newtheorem{examples}[theorem]{Examples}
\newtheorem{question}[theorem]{Question}
\newtheorem{sassumptions}[theorem]{Standing Assumptions}
\newtheorem{sassumption}[theorem]{Standing Assumption}
\newtheorem{conjecture}[theorem]{Conjecture}

\newcommand{\begintheorem}{\addtocounter{equation}{1}\begin{theorem}}
\newcommand{\beginlemma}{\addtocounter{equation}{1}\begin{lemma}}
\newcommand{\beginproposition}{\addtocounter{equation}{1}\begin{proposition}}
\newcommand{\beginsublemma}{\addtocounter{equation}{1}\begin{sublemma}}
\newcommand{\begindefinition}{\addtocounter{equation}{1}\begin{definition}}
\newcommand{\begincorollary}{\addtocounter{equation}{1}\begin{corollary}}
\newcommand{\beginproblem}{\addtocounter{equation}{1}\begin{problem}}
\newcommand{\beginremark}{\addtocounter{equation}{1}\begin{remark}}
\newcommand{\beginclaim}{\addtocounter{equation}{1}\begin{claim}}
\newcommand{\beginassumptions}{\addtocounter{equation}{1}\begin{assumptions}}
\newcommand{\beginexamples}{\addtocounter{equation}{1}\begin{examples}}
\newcommand{\beginquestion}{\addtocounter{equation}{1}\begin{question}}
\newcommand{\beginsassumptions}{\addtocounter{equation}{1}\begin{sassumptions}}
\newcommand{\beginsassumption}{\addtocounter{equation}{1}\begin{sassumption}}
\newcommand{\beginconjecture}{\addtocounter{equation}{1}\begin{conjecture}}

\begin{document}

\title{An Introduction to Heisenberg Groups}

\author{Stephen Semmes}

\date{}

\maketitle

	Actually, a number of related objects are often called
``Heisenberg groups'', and they appear in several ways.

	Let $n$ be a positive integer, and let $\lambda$ be a positive
real number.  For each $x$ in ${\bf R}^n$, consider the linear
operator $T_x$ acting on the vector space of continuous complex-valued
functions on ${\bf R}^n$ with compact support defined by
\begin{equation}
	T_x(f)(w) = f(w - x).
\end{equation}
In other words, $T_x(f)$ is the same as $f$ translated by $x$, in the
sense that $T_x(f)(w + x) = f(w)$.  Thus $T_0$ is the identity
operator, and
\begin{equation}
	T_{x + x'} = T_x \circ T_{x'}
\end{equation}
for all $x, x' \in {\bf R}^n$.  Also, $T_x$ is the identity operator
only when $x = 0$, and thus the correspondence $x \mapsto T_x$
provides an isomorphic embedding of ${\bf R}^n$, as an abelian group,
into the group of invertible linear operators on the vector space of
continuous functions with compact support on ${\bf R}^n$ under the
binary operation composition.  Of course there are numerous other
vector spaces of functions on ${\bf R}^n$ that one might use here.

	For each $y$ in ${\bf R}^n$ define a linear operator $U_y$ on
the same space of functions on ${\bf R}^n$ by
\begin{equation}
	U_y(f)(w) = \exp (2 \pi i \, \lambda \, y \cdot w) \, f(w),
\end{equation}
where $y \cdot w$ denotes the usual inner product on ${\bf R}^n$,
namely,
\begin{equation}
	y \cdot w = \sum_{j=1}^n y_j w_j, \quad y = (y_1, \ldots, y_n),
					  \ w = (w_1, \ldots, w_n).
\end{equation}
As before, $U_y$ is equal to the identity operator if and only if $y = 0$,
and
\begin{equation}
	U_{y + y'} = U_y \circ U_{y'}
\end{equation}
for all $y, y' \in {\bf R}^n$.  Thus the correspondence $y \mapsto U_y$
also defines an isomorphic embedding of ${\bf R}^n$ as an abelian group
into the group of invertible linear transformations on our vector
space of functions.  Note that the operator $U_y$ can be described
in terms of translations on the Fourier transform side, rather than
directly, as for $T_x$.
	
	How do these two families of operators interact?  What is the
subgroup of the group of invertible linear operators that they generate?
If $\alpha$ is a complex number which satisfies $|\alpha| = 1$, let
$C_\alpha$ denote the linear operator on our vector space of functions
on ${\bf R}^n$ defined by
\begin{equation}
	C_\alpha(f)(w) = \alpha \, f(w).
\end{equation}
That is, we simply multiply the function $f$ by the scalar $\alpha$.
The correspondence $\alpha \mapsto C_\alpha$ defines an isomorphic
embedding of the abelian group of complex numbers with modulus $1$,
under multiplication, into the group of invertible linear operators on
our vector space of functions.  Clearly $C_\alpha$ commutes with $T_x$
and $U_y$ for all $\alpha$, $x$, and $y$ in their respective domains,
but in general $T_x$ and $U_y$ do not commute.  However, for all $x$,
$y$ in ${\bf R}^n$, $U_y \circ T_x$ can be written as $T_x \circ U_y
\circ C_\alpha$ for a complex number $\alpha$ of modulus $1$.  Indeed,
for any function $f$ on ${\bf R}^n$,
\begin{eqnarray}
	T_x(U_y(f))(w) & = & U_y(f)(w - x) 				 \\
	  & = & \exp (2 \pi i \, \lambda \, y \cdot (w-x)) \, f(w - x)    
								\nonumber \\
	  & = & \exp (- 2 \pi i \, \lambda \, y \cdot x) \,
	     \exp (2 \pi i \, \lambda \, y \cdot w) \, T_x(f)(w)
							  	\nonumber \\
	  & = & \exp (- 2 \pi i \, \lambda \, y \cdot x) \, U_y(T_x(f))(w),
								\nonumber
\end{eqnarray}
so that
\begin{equation}
	U_y \circ T_x = T_x \circ U_y \circ C_\alpha \quad\hbox{with }
			   \alpha = \exp (2 \pi i \, \lambda \, y \cdot x).
\end{equation}

	It is not hard to check that the collection of linear
operators on our vector space of functions on ${\bf R}^n$ of the form
$T_x \circ U_y \circ C_\alpha$, where $x$, $y$ range through all
elements of ${\bf R}^n$ and $\alpha$ ranges through all complex
numbers of modulus $1$, forms a subgroup of the group of invertible
linear operators on this vector space.  This subgroup is the same as
the subgroup generated by the $T_x$'s and $U_y$'s for $x, y \in {\bf R}^n$.

	Here is a closely connected point.  For this it is convenient
to work with functions $f$ which are continuously differentiable.  If
$\nu$ is an element of ${\bf R}^n$, then we write $D_\nu$ for the
operator of directional differentiation in the direction $\nu$, i.e.,
\begin{equation}
	D_\nu(f)(w) = \frac{d}{dt} \Big|_{t = 0} \, f(w + t \, \nu).
\end{equation}
If $\mu$ is a linear function on ${\bf R}^n$, so that $\mu(x) = x \cdot u$
for some $u$ in ${\bf R}^n$ and all $x$ in ${\bf R}^n$, then let us
write $M_\mu$ for the operator of multiplication by $\mu$,
\begin{equation}
	M_\mu(f)(w) = \mu(w) \, f(w).
\end{equation}
In general $D_\nu$ and $M_\mu$ do not commute, and in fact
\begin{equation}
	D_\nu (M_\mu f) - M_\mu (D_\nu f) = \mu(\nu) \, f
\end{equation}
for all continuously differentiable functions $f$ on ${\bf R}^n$,
where the scalar $\mu(\nu)$ is the same as $D_\nu(\mu)$, which is a
constant since $\mu$ is linear.

	Let us define $H_n({\bf R})$ as follows.  First, as a set,
$H_n({\bf R})$ is equal to ${\bf R}^n \times {\bf R}^n \times {\bf
R}$, i.e., the set of ordered triples $(x, y, t)$ where $x$, $y$ lie
in ${\bf R}^n$ and $t$ lies in ${\bf R}$.  Next, we define a binary
operation on $H_n({\bf R})$ by 
\begin{equation}
	(x, y, t) \cdot (x', y', t') = (x + x', y + y', t + t' + x' \cdot y).
\end{equation}
With respect to this operation, one can verify that $H_n({\bf R})$ becomes
a group, with $(0, 0, 0)$ as the identity element, and $(-x, -y, -t)$
the inverse of $(x, y, t)$.

	Consider the correspondence
\begin{equation}
	(x, y, t) \mapsto T_x \circ U_y \circ C_\alpha, 
			\quad \alpha = \exp (2 \pi i \, \lambda \, t).
\end{equation}
This defines a mapping from $H_n(R)$ into invertible linear mappings
acting on functions on ${\bf R}^n$, and one can verify that this is
a group homomorphism.  A triple $(x, y, t)$ corresponds to the identity
operator if and only if $x = y = 0$ and $t$ is an integer multiple
of $1/\lambda$.

	We can define $H_n({\bf Z})$ to be the set of triples $(x, y, t)$
where $x$, $y$ lie in ${\bf Z}^n$ and $t$ lies in ${\bf Z}$.  It is easy
to see that this defines a subgroup of $H_n({\bf R})$.

	In addition to being a group, $H_n({\bf R})$ is a smooth manifold
of dimension $2 n + 1$, and the group operation is smooth.  The quotient
$H_n({\bf R})/H_n({\bf Z})$ makes sense not only as a set of cosets, but
also as a compact smooth manifold without boundary of dimension $2 n + 1$.
When $n = 1$, this manifold has dimension $3$, and is one of the basic
building blocks discussed in \cite{Thurston}.

	Define $a_j$ in $H_n({\bf Z})$ for $j = 1, \ldots, n$ to be
the triple $(x, y, t)$ such that all components of $x$ are equal to $0$
except for the $j$th component, which is equal to $1$, and such that
$y = 0$ and $t = 0$.  Similarly, define $b_j$ in $H_n({\bf Z})$ for 
$j = 1, \ldots, n$ to be the triple $(x, y, t)$ such that $x = 0$,
all components of $y$ are equal to $0$ except for the $j$th component,
which is equal to $1$, and $t = 0$.  Define $c$ in $H_n({\bf Z})$ to
be the triple $(0, 0, 1)$.  It is easy to see that
\begin{equation}
	a_1, \ldots, a_n, b_1, \ldots, b_n, c
\end{equation}
generate $H_n({\bf Z})$.  Indeed, if $k$, $l$ lie in ${\bf Z}^n$ and
$m$ lies in ${\bf Z}$, then $(k, l, m)$ is the same as
\begin{equation}
\label{normal form for elements of H_n({bf Z}) in terms of the generators}
	a_1^{k_1} \cdots a_n^{k_n} b_1^{l_1} \cdots b_n^{l_n} c^m
\end{equation}
in the group, where of course this expression is interpreted using the
group operation. 

	For $i, j = 1, \ldots, n$ we have that
\begin{equation}
	a_i \, a_j = a_j \, a_i, \quad b_i \, b_j = b_j \, b_i
\end{equation}
and
\begin{equation}
	a_i \, c = c \, a_i, \quad b_j \, c = c \, b_j.
\end{equation}
If $1 \le i, j \le n$ and $i \ne j$, then
\begin{equation}
	a_i \, b_j = b_j \, a_i.
\end{equation}
When $i = j$ we have in place of this
\begin{equation}
	b_i \, a_i = a_i \, b_i \, c
\end{equation}
for $i = 1, \ldots, n$.  These relations describe the group $H_n({\bf
Z})$ completely; every element of $H_n({\bf Z})$ can be represented in
a unique way as (\ref{normal form for elements of H_n({bf Z}) in terms
of the generators}), and these relations are adequate to define the
group operation in terms of this representation.

	In addition to the group structure on $H_n({\bf R})$, there is
a natural family of \emph{dilations}.  Namely, for each positive real
number $r$, define a mapping $\delta_r$ from $H_n({\bf R})$ to itself
by
\begin{equation}
	\delta_r(x, y, t) = (r \, x, r \, y, r^2 \, t).
\end{equation}
This is a one-to-one mapping of $H_n({\bf R})$ onto itself, $\delta_r
\circ \delta_s = \delta_{r \, s}$ for all $r, s > 0$, and $\delta_1$
is the identity mapping on $H_n({\bf R})$.  It is easy to check that
these dilations preserve the group structure on $H_n({\bf R})$.

	In connection with complex variables it is convenient to
employ a slightly different formulation.  Define $\widetilde{H}_n$ as
a set as ${\bf C}^n \times {\bf R}$.  The group operation on
$\widetilde{H}_n$ is defined by
\begin{equation}
	(z, t) \cdot (z', t') 
   = \biggl(z + z', t + t' + 2 \im \sum_{j=1}^n z_j \, \overline{z'_j} \biggr),
\end{equation}
where $z = (z_1, \ldots, z_n)$, $z' = (z'_1, \ldots, z'_n)$, and
$\overline{a}$ denotes the complex conjugate of a complex number $a$.
As before, $(0,0)$ is the identity element of $\widetilde{H}_n$, and
$(-z, -t)$ is the inverse of $(z, t)$.

	Put
\begin{equation}
	U = \{(w, \sigma) \in {\bf C}^n \times {\bf C} : \im \sigma > |w|^2\},
\end{equation}
where $|w|^2 = \sum_{j=1}^n |w_j|^2$, as usual.  Thus
\begin{equation}
	\partial U = \{(w, \sigma) \in {\bf C}^n \times {\bf C} :
							\im \sigma = |w|^2\}.
\end{equation}
For each $(z,t)$ in $\widetilde{H}_n$, define the mapping $A_{(z, t)}$
from $\overline{U}$ to itself by
\begin{equation}
	A_{(z, t)}(w,\sigma) = 
   \biggl(w + z, \sigma + t + i |z|^2 + 
			2 i \sum_{j=1}^n w_j \, \overline{z_j} \biggr).
\end{equation}
Because
\begin{equation}
	|w + z|^2 = |w|^2 + 2 \re \sum_{j=1}^n w_j \, \overline{z_j} + |z|^2,
\end{equation}
one can check that $A_{(z, t)}(w, \sigma) \in U$ when $(w, \sigma) \in
U$ and $A_{(z, t)}(w, \sigma) \in \partial U$ when $(w, \sigma) \in
\partial U$.  Also,
\begin{eqnarray}
\lefteqn{\qquad A_{(z, t)} (A_{(z', t')}(w, \sigma))} && 		\\
	& = & \biggl(w + z' + z, \sigma + t' + t + i |z'|^2 + i |z|^2
			+ 2 i \sum_{j=1}^n w_j \, \overline{z'_j}
		+ 2 i \sum_{j=1}^n (w_j +z'_j) \, \overline{z_j} \biggr)
								\nonumber \\
	& = & \biggl(w + z' + z, \sigma + t' + t + i |z'|^2 + i |z|^2
			+ 2 i \sum_{j=1}^n w_j \, \overline{(z'_j + z_j)}
		+ 2 i \sum_{j=1}^n z'_j \, \overline{z_j} \biggr)
								\nonumber \\
	& = & \biggl(w + z' + z, \sigma + t' + t + i |z' + z|^2
			+ 2 i \sum_{j=1}^n w_j \, \overline{(z'_j + z_j)}
		- 2  \im \sum_{j=1}^n z'_j \, \overline{z_j} \biggr)
								\nonumber \\
	& = & \biggl(w + z + z', \sigma + t + t' + i |z + z'|^2
			+ 2 i \sum_{j=1}^n w_j \, \overline{(z_j + z'_j)}
		+ 2  \im \sum_{j=1}^n z_j \, \overline{z'_j} \biggr)
								\nonumber \\
	& = & A_{(z, t) \cdot (z', t')}(w, \sigma).		\nonumber
\end{eqnarray}
In particular, $A_{(-z, -t)}$ is the inverse of $A_{(z, t)}$ as a mapping
from $\overline{U}$ onto itself.

	From the definition of $A_{(z, t)}$ it is clear that it is a
\emph{holomorphic} mapping of $U$ to itself, and in fact a
complex-affine mapping, since there are only $w$'s and $\sigma$'s in
$A_{(z,t)}(w, \sigma)$, and none of their complex conjugates.  Thus $A_{(z,
t)}$ is in fact a \emph{biholomorphic} mapping of $U$ onto itself.

	Let us abuse our notation a bit and use $\delta_r$ to denote
the natural dilation by $r > 0$ on $\widetilde{H}_n$,
\begin{equation}
	\delta_r(z, t) = (r \, z, r^2 \, t).
\end{equation}
As before, $\delta_r$ is a one-to-one mapping of $\widetilde{H}_n$
onto itself, $\delta_r \circ \delta_s = \delta_{r \, s}$ for all $r, s
> 0$, $\delta_1$ is the identity mapping on $\widetilde{H}_n$, and
$\delta_r$ preserves the group structure on $\widetilde{H}_n$.  These
dilations also correspond to mappings on $\overline{U}$ in a natural
way.  Namely, define $\Delta_r : \overline{U} \to \overline{U}$ for $r
> 0$ by
\begin{equation}
	\Delta_r(w, \sigma) = (r \, w, r^2 \, \sigma).
\end{equation}
It is easy to see that $\Delta_r$ is a one-to-one mapping of $\overline{U}$
onto itself which takes $U$ to $U$, and that $\Delta_r \circ \Delta_s
= \Delta_{r \, s}$ for all $r, s > 0$ and $\Delta_1$ is the identity
mapping.  Furthermore, $\Delta_r$ is a biholomorphic mapping of $U$ onto
itself, and
\begin{equation}
	\Delta_r(A_{(z, t)}(w, \sigma)) 
		= A_{\delta_r(z, t)}(\Delta_r(w, \sigma)).
\end{equation}

	A famous fact about $U$ is that it is biholomorphically
equivalent to the unit ball in ${\bf C}^{n+1}$.  Thus the $A_{(z,
t)}$'s and $\Delta_r$'s correspond to biholomorphic transformations on
the unit ball in ${\bf C}^{n+1}$.  They do not account for all of the
biholomorphic transformations, though; basically what is missing are
the complex-linear tranformations on ${\bf C}^{n+1}$ which map the
unit ball onto itself, which is the same as saying that they preserve
the standard Hermitian inner product on ${\bf C}^{n+1}$.


\begin{thebibliography}{99}


\bibitem {Alexopoulos} G.~Alexopoulos, {\it Sub-Laplacians with
Drift on Lie Groups of Polynomial Volume Growth}, Memoirs of the
American Mathematical Society {\bf 739}, 2002.

\bibitem {A-T} L.~Ambrosio and P.~Tilli, {\it Selected Topics on
``Analysis on Metric Spaces''}, Scuola Normale Superiore, Pisa, 2000.

\bibitem {B-Ga-Gr} R.~Beals, B.~Gaveau, and P.~Greiner, {\it
Hamilton--Jacobi theory and the heat kernel on Heisenberg groups},
Journal de Math\'ematiques Pures et Appliqu\'ees (9) {\bf 79} (2000),
633-689.

\bibitem {B-Gr} R.~Beals and P.~Greiner, {\it Calculus on Heisenberg
Manifolds}, Annals of Mathematics Studies {\bf 119}, Princeton University
Press, 1988.

\bibitem {B-Gr-S1} R.~Beals, P.~Greiner, and N.~Stanton, {\it 
The heat equation and geometry of CR manifolds}, Bulletin of the American
Mathematical Society (N.S.) {\bf 10} (1984), 275--276.

\bibitem {B-Gr-S2} R.~Beals, P.~Greiner, and N.~Stanton, {\it
The heat equation on a CR manifold}, Journal of Differential Geometry
{\bf 20} (1984), 343--387.

\bibitem {B-R} A.~Bella\"{\i}che and J.-J.~Risler, editors, {\it
Sub-Riemannian Geometry}, Birkh\"auser, 1996.

\bibitem {B-C-T} C.~Berenstein, D.-C.~Chang, and J.~Tie, {\it Laguerre
Calculus and its Applications on the Heisenberg Group}, American
Mathematical Society and International Press, 2001.

\bibitem {Borel} A.~Borel, {\it Semisimple Groups and Riemannian
Symmetric Spaces}, Hindustan Book Agency, 1998.

\bibitem {Burillo} J.~Burillo, {\it Lower bounds of isoperimetric
functions for nilpotent groups}, in {\it Geometric and Computational
Perspectives on Infinite Groups}, DIMACS Series in Discrete
Mathematics and Theoretical Computer Science {\bf 25}, American
Mathematical Society, 1996.

\bibitem {C-W1} R.~Coifman and G.~Weiss, {\it Analyse Harmonique
Non-Commutative sur Certains Espaces Homog\`enes}, Lecture Notes in
Mathematics {\bf 242}, Springer-Verlag, 1971.

\bibitem {C-W2} R.~Coifman and G.~Weiss, {\it Extensions of Hardy
spaces and their use in analysis}, Bulletin of the American Mathematical
Society {\bf 83} (1977), 569--645.

\bibitem {C-SC} T.~Coulhon and L.~Saloff-Coste, {\it
Isop\'erim\'etrie pours les groupes et les vari\'et\'es}, Revista
Matem\'atica Iberoamericana {\bf 9} (1993), 293--314.

\bibitem {D-G-N} D.~Danielli, N.~Garofalo, and D.-M.~Nhieu, {\it Notions
of convexity in Carnot groups}, preprint, 2002.

\bibitem {ECHLPT} D.~Epstein, J.~Cannon, D.~Holt, S.~Levy,
M.~Paterson, and W.~Thurston, {\it Word Processing in Groups}, Jones and
Bartlett, 1992.

\bibitem {Benson} B.~Farb, {\it Automatic groups: A guided tour},
L'Enseignement Math\'ematique (2) {\bf 38} (1992), 291--313.

\bibitem {Folland} G.~Folland, {\it Harmonic Analysis in Phase Space},
Annals of Mathematics Studies {\bf 122}, Princeton University Press,
1989.

\bibitem {F-S} G.~Folland and E.~Stein, {\it Hardy Spaces on
Homogeneous Groups}, Mathematical Notes {\bf 28}, Princeton University
Press, 1982.

\bibitem {F-S-SC1} B.~Franchi, R.~Serapioni, and F.~Serra Cassano,
{\it Sur les ensembles de p\'erim\`etre fini dans le groupe de Heisenberg},
Comptes Rendus des S\'eances de l'Acad\'emie des Sciences Paris S\'erie I
Math\'ematiques {\bf 329} (1999), 183--188.

\bibitem {F-S-SC2} B.~Franchi, R.~Serapioni, and F.~Serra Cassano,
{\it Rectifiability and perimeter in the Heisenberg group}, Mathematische
Annalen {\bf 321} (2001), 479--531.

\bibitem {F-S-SC3} B.~Franchi, R.~Serapioni, and F.~Serra Cassano,
{\it Regular hypersurfaces, intrinsic perimeter and implicit function
theorem in Carnot groups}, to appear, Communications in Analysis and
Geometry.

\bibitem {Goldman} W.~Goldman, {\it Complex Hyperbolic Geometry},
Oxford University Press, 1999.

\bibitem {H-K} P.~Hajlasz and P.~Koskela, {\it Sobolev Met Poincar\'e},
Memoirs of the American Mathematical Society {\bf 688}, 2000.

\bibitem {dlH} P.~de la Harpe, {\it Topics in Geometric Group Theory},
University of Chicago Press, 2000.

\bibitem {Juha1} J.~Heinonen, {\it Calculus on Carnot groups}, in {\it
Fall School in Analysis (Jyv\"askyl\"a, 1994)}, 1--31, Reports of the
Department of Mathematics, University of Jyv\"askyl\"a {\bf 68}, 1995.

\bibitem {Juha2} J.~Heinonen, {\it Lectures on Analysis on Metric Spaces},
Springer-Verlag, 2001.

\bibitem {Helgason} S.~Helgason, {\it Differential Geometry, Lie Groups,
and Symmetric Spaces}, Academic Press, 1978.

\bibitem {Howe} R.~Howe, {\it On the role of the Heisenberg group in
harmonic analysis}, Bulletin of the American Mathematical Society (N.S.)
{\bf 3} (1980), 821--843.

\bibitem {Kapovich} M.~Kapovich, {\it Hyperbolic Manifolds and
Discrete Groups}, Birkh\"auser, 2001.

\bibitem {Kobayashi} S.~Kobayashi, {\it Transformation Groups in
Differential Geometry}, Springer-Verlag, 1972.

\bibitem {K-N} S.~Kobayashi and K.~Nomizu, {\it Foundations of
Differential Geometry, Volume II}, Wiley, 1969.

\bibitem {Koranyi1} A.~Koranyi, {\it The Poisson integral for
generalized half-planes and bounded symmetric domains}, Annals of
Mathematics (2) {\bf 82} (1965), 332--350.

\bibitem {Koranyi2} A.~Koranyi, {\it Harmonic functions on
Hermitian hyperbolic space}, Transactions of the American Mathematical
Society {\bf 135} (1969), 507--516.

\bibitem {Koranyi3} A.~Koranyi, {\it A remark on boundary values
of functions of several variables}, in {\it Several Complex Variables I,
Maryland 1970}, Lecture Notes in Mathematics {\bf 155}, Springer-Verlag, 
1970.

\bibitem {K-R1} A.~Koranyi and H.~Reimann, {\it Quasiconformal
mappings on the Heisenberg group}, Inventiones Mathematicae {\bf 80}
(1985), 309--338.

\bibitem {K-R2} A.~Koranyi and H.~Reimann, {\it Foundation for
the theory of quasiconformal mappings on the Heisenberg group},
Advances in Mathematics {\bf 111} (1995), 1--87.

\bibitem {K-V1} A.~Koranyi and S.~Vagi, {\it Int\'egrales
singuli\`eres sur certains espaces homog\`enes}, Comptes Rendus
Hebdomadiares des S\'eances de l'Acad\'emie des Sciences S\'eries A et
B {\bf 268} (1969), A765--A768.

\bibitem {K-V2} A.~Koranyi and S.~Vagi, {\it Singular integrals
in homogeneous spaces and some problems in classical analysis}, Annali
della Scuola Normale Superiore di Pisa (3) {\bf 25} (1971), 575--648.

\bibitem {K-V3} A.~Koranyi and S.~Vagi, {\it Cauchy--Szeg\"o
integrals for systems of harmonic functions}, Annali
della Scuola Normale Superiore di Pisa (3) {\bf 26} (1972), 181--196.

\bibitem {K-V-W} A.~Koranyi, S.~Vagi, and G.~Welland, {\it Remarks
on the Cauchy integral and the conjugate function in generalized half-planes},
Journal of Mathematics and Mechanics {\bf 19} (1970) 1069--1081.

\bibitem {Krantz1} S.~Krantz, {\it Function Theory of Several
Complex Variables}, Wiley, 1982.

\bibitem {Krantz2} S.~Krantz, {\it Complex Analysis: The
Geometric Viewpoint}, Carus Mathematical Monographs {\bf 23},
Mathematical Association of America, 1990.

\bibitem {Krantz3} S.~Krantz, {\it A Panorama of Harmonic Analysis},
Carus Mathematical Monographs {\bf 27}, Mathematical Association of
America, 1999.

\bibitem {L-R} G.-P.~Leonardi and S.~Rigot, {\it Isoperimetric sets on
Carnot groups}, to appear, Houston Journal of Mathematics.

\bibitem {L-M-S} G.~Lu, J.~Manfredi, and B.~Stroffolini, {\it Convex
functions on the Heisenberg group}, preprint, 2002.

\bibitem {Magnani} V.~Magnani, {\it Elements of Geometric Measure
Theory on Sub-Riemannian Groups}, thesis, Scuola Normale Superiore,
Pisa, 2002.

\bibitem {Mitchell} J.~Mitchell, {\it On Carnot--Carath\'eodory
metrics}, Journal of Differential Geometry {\bf 21} (1985), 35--45.

\bibitem {N-S} A.~Nagel and E.~Stein, {\it Lectures on Pseudodifferential
Operators: Regularity Theorems and Applications to Nonelliptic Problems},
Mathematical Notes {\bf 24}, Princeton University Press, 1979.

\bibitem {Pierre1} P.~Pansu, {\it Une in\'egalit\'e isop\'erim\'etrique
sur le groupe de Heisenberg}, Comptes Rendus des S\'eances de l'Acad\'emie
des Sciences S\'erie I Math\'ematiques {\bf 295} (1982), 127--130.

\bibitem {Pierre2} P.~Pansu, {\it An isoperimetric inequality on
the Heisenberg group}, Conference on Differential Geometry on Homogeneous
spaces (Torino, 1983), Rendiconti del Seminario Matematico, Universit\`a
e Politecnico di Torino, Special Issue (1983), 159--174.

\bibitem {Pierre3} P.~Pansu, {\it M\'etriques de
Carnot--Carath\'eodory et quasiisom\'etries des espaces sym\'etriques
de rang un}, Annals of Mathematics (2) {\bf 129} (1989), 1--60.

\bibitem {Rud1} W.~Rudin, {\it Functional Analysis},
McGraw-Hill, 1973.

\bibitem {Rud2} W.~Rudin, {\it Function Theory in the Unit Ball
of ${\bf C}^n$}, Springer-Verlag, 1980.

\bibitem {Ste1} E.~Stein, {\it Singular Integrals and
Differentiability Properties of Functions}, Princeton University
Press, 1970.

\bibitem {Ste2} E.~Stein, {\it Boundary Behavior of Holomorphic
Functions of Several Complex Variables}, Mathematical Notes {\bf 11},
Princeton University Press, 1972.

\bibitem {Ste3} E.~Stein, {\it Harmonic Analysis: Real-Variable
Methods, Orthogonality, and Oscillatory Integrals}, Princeton
University Press, 1993.

\bibitem {S-W} E.~Stein and G.~Weiss, {\it Introduction to
Fourier Analysis on Euclidean Spaces}, Princeton University Press,
1971.

\bibitem {Strichartz} R.~Strichartz, {\it Sub-Riemannian
Geometry}, Journal of Differential Geometry {\bf 24} (1986), 221--263;
corrections, {\bf 30} (1989), 595--596.

\bibitem {Taylor1} M.~Taylor, {\it Noncommutative microlocal analysis I},
Memoirs of the American Mathematical Society {\bf 313}, 1984.

\bibitem {Taylor2} M.~Taylor, {\it Noncommutative Harmonic Analysis},
Mathematical Surveys and Monographs {\bf 22}, American Mathematical
Society, 1986.

\bibitem {Thurston} W.~Thurston, {\it Three-Dimensional Geometry
and Topology}, edited by S.~Levy, Princeton University Press, 1997.

\bibitem {Vagi} S.~Vagi, {\it Harmonic analysis on Cartan and
Siegel domains}, it {\it Studies in Harmonic Analysis}, MAA Studies in
Mathematics {\bf 13}, 257--309, Mathematical Association of America,
1976.

\bibitem {V-M} N.~Varopoulos and S.~Mustapha, {\it Analysis on
Lie Groups}, Cambridge University Press, 2001.

\bibitem {V-SC-C} N.~Varopoulos, L.~Saloff-Coste, and T.~Coulhon,
{\it Analysis and Geometry on Groups}, Cambridge University Press, 1992.



\end{thebibliography}
\end{document}